\documentclass[preprint,12pt]{elsarticle}
\usepackage{multirow}
 
\usepackage{graphicx}        
\usepackage{multicol}        
\usepackage[bottom]{footmisc}
\usepackage{natbib}
\usepackage{url}

\usepackage{epsfig}

\renewcommand{\PackageWarningNoLine}[2]{}
\usepackage{amsmath}
\usepackage{amsfonts}

\biboptions{sort&compress}
\journal{Journal of Computational and Applied Mathematics}

\begin{document}

 \begin{frontmatter}

\title{Numerical upscaling of the  free boundary dam problem in multiscale high-contrast media}

\author{\textbf{Juan Galvis}$^{1*}$}
\cortext[cor1]{Email address : jcgalvisa@unal.edu.co}

\author{\textbf{Luis F. Contreras}$^{1,2**}$}
\cortext[cor1]{Email address : lfcontrerash@unal.edu.co, luis.contreras@escuelaing.edu.co}

\author{\textbf{Carlos V\'azquez}$^{3}$}
\cortext[cor1]{Email address : carlosv@udc.es}

\address{$^{1}$ Departamento de Matem\'{a}ticas \\
Universidad Nacional de Colombia \\ 
Bogot\'{a} D.C.}

\address{$^{2}$ Departamento de Matem\'{a}ticas \\
Escuela C Julio Garavito.\\
Bogot\'{a} D.C. }

\address{$^{3}$ Departamento de Matem\'{a}ticas \\
Universidade da Coruna \\ 
A Coruna}

\begin{abstract}
In this paper, we address the numerical homogenization approximation of a free-boundary dam problem posed in a heterogeneous media. More precisely, we propose a generalized multiscale finite element (GMsFEM) method for the heterogeneous dam problem. The motivation of using the GMsFEM approach comes from the multiscale nature of the porous media due to its high-contrast permeability. Thus, although we can classically formulate the free-boundary dam problem as in the homogeneous case, a very high resolution will be needed by a standard finite element approximation in order to obtain realistic results that recover the multiscale nature. First, we introduce a fictitious time variable which motivates a suitable time discretization that can be understood as a fixed point iteration to the steady state solution, and we use a duality method to deal with the involved multivalued nonlinear terms. Next, we compute  efficient approximations of the pressure and the saturation by using the GMfsFEM method and we can identify the free boundary. More precisely, the GMsGEM method provides numerical results that capture the behavior of the solution due to the variations of the coefficient at the fine-resolution, by just solving linear systems with size proportional to the number of coarse blocks of a  coarse-grid (that does not need to be adapted to the variations of the coefficient). Finally, we present illustrative numerical results to validate the proposed methodology. 
 \end{abstract}

 \end{frontmatter}
 
\section{Introduction}
In this paper, we start from the heterogeneous dam problem  originally posed in  \cite{MR3062565}. More precisely, in \cite{MR3062565}  the numerical homogenization approximation of a free-boundary dam problem posed in a heterogeneous media with scale separation is considered. In the here presented work, we deal with the numerical upscaling of a similar free-boundary problem posed on a high-contrast multiscale media, in this case with no-scale separation assumption. Thus, although we can formulate the free-boundary problem, due to the multiscale nature of the porous media, a very high resolution will be needed by a finite element approximation in order to obtain realistic results.

Following \cite{BermudezDurany, MR3062565}, we first approximate the nonlinear steady state dam problem by a time dependent one coming form the introduction of a fictitious time variable, this procedure can be understood as a fixed point iteration indexed by the fictitious time variable. For the time discretization, we consider a characteristics method, which is based on the numerical approximation of material (or total) derivative, that is a concept well understood in continuum mechanics. Moreover, a duality method is considered to deal with the involved multivalued nonlinear terms. As in \cite{BermudezDurany, MR3062565}, after applying these techniques, at each iteration the spatial approximation of a resulting pressure equation is required. 

The main innovative achievement of the present work comes from the proposed method to address the spatial approximation of the pressure in the case of no-scale separation in the heterogeneous porous medium. In fact, solving for the pressure equation at the same resolution of the medium is described results to be impractical for this application, as it is the case in several porous media flow models involving multiple scales. Here, we originally propose to compute an efficient approximation of the pressure by  employing the   Generalized Multiscale Finite Element approximation introduced \cite{egh12,eglp13oversampling,MR3477310} and references therein. The GMsFEM method provides numerical results that capture the behavior  of the solution due to the variations of the coefficient at the fine-resolutions, by just solving linear systems with size proportional to the number of coarse blocks of a  coarse-grid (that does not need to be adapted to the variations of the coefficient).
 
The GMsFEM methodology has been successfully employed for the numerical upscaling and/or preconditioning of problems with a complicated dependence of physical parameters such as high-contrast and  variations over several scales. We mention some applications on heterogeneous high-contrast media 
 such as multiphase porous media flow \cite{MR3430146}, Brinkman flow \cite{MR3296062},  fracture porous media flow 
 \cite{MR3544033,MR3723657}, wave propagation, non-linear parabolic problems \cite{MR3261091},
history matching problems, nonlinear elliptic problems, flow on perforated domains, among others. 
 
Upscaling the free boundary model for the porous dam problem is a very challenging task. We consider the formulation of heterogeneous dam problem presented in \cite{BrezisKinderlehrerStampacchia, Alt, Rodrigues, Baiocchi}. In some of these references, a homogenized problem have been derived for some families of permeabilities with scale separation. We mention, in particular, \cite{MR3062565} where a homogenized problem for 
isotropic permeability coefficients $\kappa(x,x/\epsilon)$ depending on a small parameter $\epsilon$ is derived and numerically verified. Analogous problems involving pressure and saturation giving rise to homogenized models appear in the domain of lubrication of rough surfaces with cavitation phenomena, see \cite{BMV05a, BMV05b, BMV05c} and the references therein, for example. For numerical computations in the dam problem with homogeneous media, we can refer to \cite{BermudezDurany}, among other works. 
 
In some practical cases we do not have coefficients with scale separations, so that writing a homogenized model as in \cite{MR3062565} results to be an impossible task. We can still try to approximate solutions using numerical upscaling techniques that use two (or more) grids. One fine grid where all the scales and variations of the coefficients are resolved (but where computing becomes unpractical) and other coarse grid where  practical computations can be carried out, although that usually does not resolve the variations of the coefficients. One 
way to approximate solutions is to use numerical upscaling techniques that project solutions on coarse-subspaces generated by specially designed coarse-basis functions. In particular, we mention the multiscale finite element method that uses one coarse basis functions per coarse node, see \cite{eh09}. 

The multiscale finite element method is equivalent to numerical homogenization when both techniques are valid approximations. However, instead of computing effective properties of the medium, it computes multiscale basis functions that capture the behavior of a reference solution (the fine-grid solution).
Recall that in the case considered in this paper the coefficient $k$ has high-variation and discontinuities (not necessarily aligned with the coarse grid). For this problem, it is known that a higher order approximation is needed in the sense that if we use one-coarse basis functions per node (or one effective coefficient per coarse block) the resulting approximation is poor and cannot be used in practical application. Indeed,
in some cases, robust approximation properties which are independent of the contrast are required. For instance, see \cite{ge09_1, ge09_1reduceddim, Efendiev_GKiL_12} where it is demonstrated that classical numerical upscaling methods (\cite{eh09}) 
do not render robust approximation properties in terms of the contrast and multiscale variations (when no scale separation is considered). Furthermore, it is shown that one basis functions per coarse node (with the usual support) is not enough to construct adequate coarse spaces \cite{ge09_1reduceddim, MR2861243}. 

The GMsFEMs methodology aims to construct  coarse spaces for Multiscale Finite Element Methods (MsFEMs) that result in accurate coarse-scale solutions for the case of high-contrast multiscale problems and, in general, for problems with a dependence on a physical parameter that negatively affects the performance of classical numerical methods. This methodology was first developed in \cite{egh12, eglp13oversampling} 
based on some previous works \cite{ge09_1, ge09_1reduceddim, Efendiev_GKiL_12, EGG_MultiscaleMOR, Review}. 

A main ingredient in the construction is the use of an approximation of  local eigenvectors (of carefully selected local eigenvalue problems)  to construct the coarse spaces. Instead of using one coarse function per coarse node as in classical MsFEM, in the GMsFEM it was proposed to use several multiscale basis functions per coarse node. These basis functions represent important features of the solution within a coarse-grid block and they are computed using eigenvectors of a local eigenvalue problem.

In the present work, we show that the GMsFEM method can be used to numerically approximate the free boundary of the heterogeneous multiscale dam problem. Our reference solution is computed by using the numerical procedure introduced in \cite{MR3062565}. This numerical scheme is an adaption   of the numerical techniques proposed in \cite{BermudezDurany} that in turns is based on the application of characteristics methods to steady state convection-diffusion equations
with a nonlinear convection term (see \cite{BermudezDurany,MR3062565}  and references therein).

The rest of the paper is organized as follows. In Section \ref{sec:problem} we detail the problem formulation and the time discretization of the artificial auxiliar time-dependent problem. In Sections \ref{sec:evolutive} and \ref{sec:duality}  we review the introduction of the duality method for the nonlinear terms. In Section \ref{gmsfem} we present the GMsFEM method as applied to the heterogeneous dam problem. Finally, in Section \ref{sec:numerics} we present numerical  evidence of the good performance of the GMsFEM methodology for the free boundary dam problem in heterogeneous multiscale media. 


\section{Problem formulation and some numerical methods}\label{sec:problem}
In this section we first pose the heterogeneous dam problem.  Moreover, we introduce an auxiliar and artificial time dependent problem and its time discretization. Also duality methods for solving the nonlinear terms are described. In both techniques, we follow the ideas in \cite{BermudezDurany, MR3062565}.

\subsection{A free boundary dam problem in high-contrast multiscale media}
\label{sec:2.1}
\begin{figure}[h]
\centering
	\includegraphics[width = .5\textwidth, keepaspectratio = true]{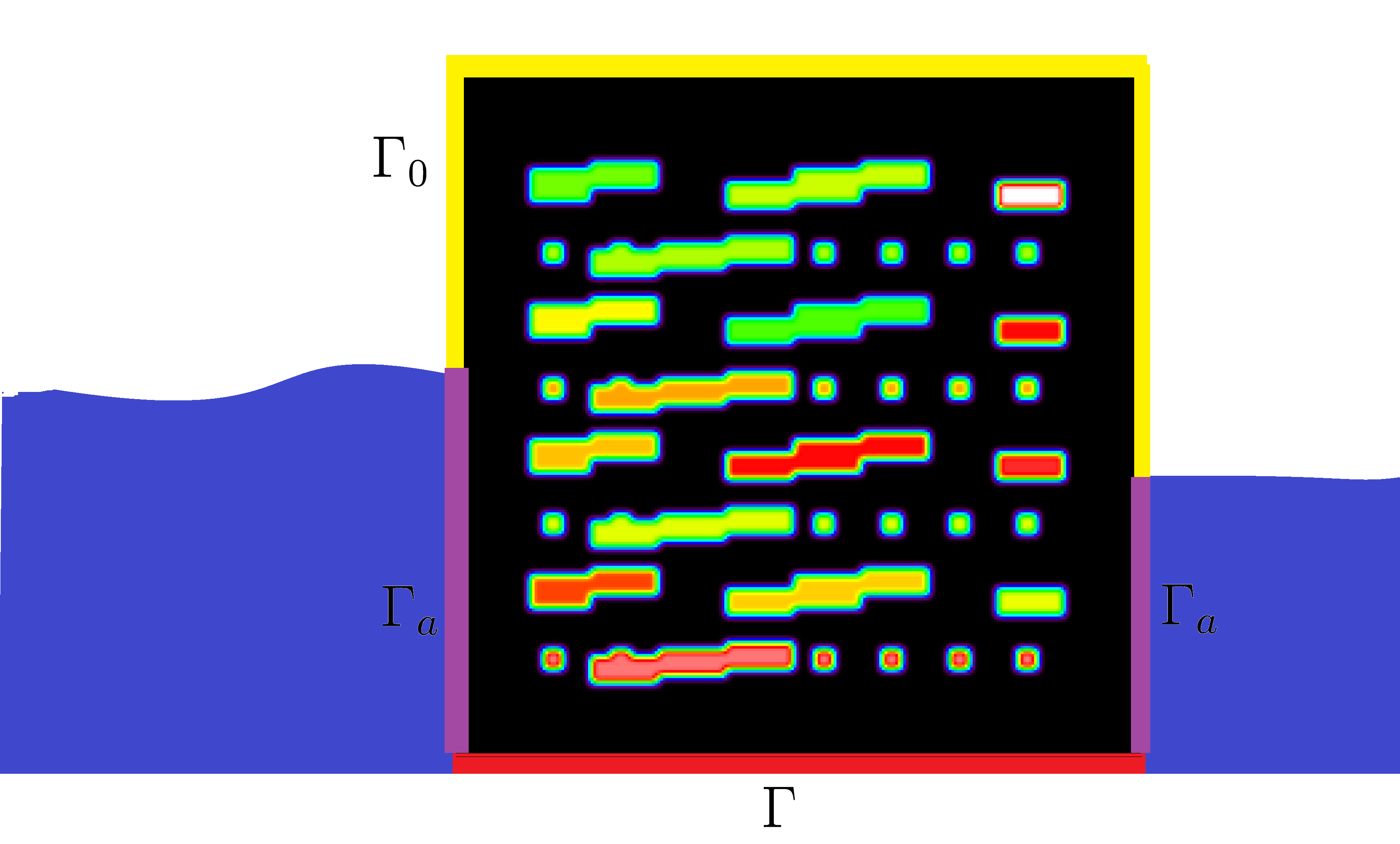}\caption{Illustration of 
	a free boundary dam problem in multiscale high-contrast porous media,}\label{illustration}
\end{figure}

In order to pose the dam problem, we consider a bounded two dimensional rectangular domain $D$ and let  $\partial D=\Gamma\cup \Gamma_0\cup \Gamma_a$ denote its boundary, where $\Gamma$  is an impervious part of the boundary, $\Gamma_0$ is the part of the 
boundary in contact with open air,  and  $\Gamma_a$ is the part of the boundary in contact with water. See Figure \ref{illustration} 
for an illustration. In the proposed dam problem, we aim to compute the pressure $p$  and the saturation $\theta$ of water, both defined on $D$, as well as to identify the free boundary separating the saturated and non saturated regions of the dam. Moreover, we denote by $\kappa$ the functional coefficient that represents the permeability of the porous media and let $\boldsymbol{g}:=-g \boldsymbol{e}_2$  denote the gravity. By using Darcy's law for porous media and the relation between pressure and water saturation, we obtain 
\begin{equation}\label{steady_dam}
-g\partial_2(\theta\kappa) -\mbox{div}(\kappa\nabla p) =0,  \quad p\geq 0, \quad \theta \in H(p), 
\end{equation}
where $H(\cdot)$ denotes the multivalued Heaviside operator, so that for positive pressure ($p>0$) the porous media is fully saturated ($\theta=1$) and $\theta \in [0,1)$ when $p=0$ in the non saturated region. In order to pose the strong formulation of the free-boundary dam problem, the set of equations (\ref{steady_dam}) is completed with the following boundary conditions: 
\begin{itemize}
	\item $p=h_a-x_2$  on $\Gamma_a$,  with $h_a$  height of the water level in contact with $\Gamma_a$,
	\item   $p=0$ in $\Gamma_0$,
	\item $\big( \theta \kappa \boldsymbol{g} -\kappa\nabla  p \big)\cdot \pmb{n}\geq 0$ in $\Gamma_0$ where we recall that $\boldsymbol{g}=-g\boldsymbol{e}_2$,
	\item $\Big( \theta \kappa \boldsymbol{g} -\kappa\nabla p\Big)\cdot \pmb{n}=0$ in $\Gamma$,
\end{itemize}
in previous equations $\pmb{n}$ represents the unitary outwards normal vector to the boundary $\partial D$.  See Figure \ref{illustration} for an illustration of the domain and different boundaries.

We remark that  the methodology proposed in this paper can be applied to general coefficients (see \cite{EfendievGLWESAIM12,egh12} and related works). 

Now, we focus on the case of 
 high-contrast multiscale coefficients. 
More precisely, we consider piece-wise smooth permeabilities.
We assume that the domain $D$ is the union of finitely many sub-domains, that is, 
\[ 
D= \bigcup_{i=1}^{N_{S}} D_i
\]
where $\{ D_i\}$ is a non-overlapping decomposition of $D$. The permeability coefficient can be written as, 
\begin{equation}\nonumber
\kappa(x)=k_i(x) \mbox{ for } x\in D_i, 
\end{equation}
$i=1,2,\dots,N_S,$ where $\kappa_i$ is a bounded smooth function (that may have oscillations) in $D_i$.
We say that the coefficient 
$\kappa(x)$ is a multiscale coefficient if it has oscillations at different scales in each subdomain $D_i$.  We also say that $\kappa$ is a high-contrast coefficient if the ratio $\eta_\Omega=\max_{x,y\in \Omega} \kappa(x)/\kappa(y) >> 1$ for some subsets 
$\Omega\subset D$ (in this case we say that $\Omega$ is a high-contrast sub-region).   As it is shown in \cite{ge09_1,ge09_1reduceddim, egw10, EfendievGLWESAIM12, EGG_MultiscaleMOR, Review} the up-scaling of pressures and fluxes of high-contrast multiscale coefficients is a challenging task where the complications are due to the local oscillations (that will generate oscillations in pressures and fluxes) and to the 
high-contrast (that will generate high velocities inside regions where permeability is high and also almost constant pressures in these regions); see \cite{egw10} for a more complete explanation. In particular, a difficult case arises when the coarse mesh (used in the numerical upscaling) is not aligned with the discontinuities of the coefficients; 
see \cite{ge09_1} and references therein. Practical situations are in flow problems where the coefficient representing the permeability involves long shaped sub regions (or channels) where the permeability is high.

In this article we consider precisely the free boundary dam problem for this case: a high-contrast multiscale coefficient where the upscaling coarse mesh is not aligned with the the discontinuities or with the oscillations of the coefficient. Therefore, apart from the complications described above, we deal with the non-linearities of the free boundary problem.

\subsection{Auxiliar evolutive problem and time discretization}\label{sec:evolutive}
The numerical computation of solutions to the previously described dam problem has been performed by several methods. For homogeneous porous medium, Alt \cite{Alt} and Marini \& Pietra \cite{MP86} built a numerical method based upon a discrete analogue of the continuous problem, by introducing finite-dimensional spaces and a fixed-point procedure. Berm\'udez \& Durany \cite{BermudezDurany} proposed to solve a transient version of the dam problem, using a combination of the method of characteristics and the finite element method~; then the solution of the nonlinear discretized problem is obtained by using a duality iterative algorithm. This approach has been adapted in \cite{MR3062565} in the case of heterogeneous porous medium for solving the resulting homogenized problems and comparing its solution with the small parameter dependent problem associated to the stratified medium.

In the present article we follow the methodology proposed in  \cite{BermudezDurany} and extended in \cite{MR3062565}. For this purpose, we introduce an artificial dependence on time in all the involved unknowns, so that $p=p(t,x)$ and $\theta=\theta(t,x)$, which we will denote as in the steady case with a certain abuse of notation. Therefore, we write (\ref{steady_dam}) equivalently in terms of the following evolution equations: 
\begin{equation}\label{evolutive_dam}
\frac{\partial}{\partial t}(\theta \kappa)-g\partial_2(\theta\kappa) -\mbox{div}(\kappa\nabla p) =0,  \quad p\geq 0, \quad \theta \in H(p). 
\end{equation}
Note that the first equation in (\ref{evolutive_dam}) is a nonlinear advection-diffusion equation. In order to discretize in time this equation by using the characteristics method, we first introduce the material or total derivative associated to the vector field driving the convection term which is given by $\boldsymbol{g}=(0,-g)$, so that  
$$\displaystyle \frac{D z}{D t} =\frac{\partial z}{\partial t} +{\mathbf g} \cdot \nabla z=\frac{\partial z}{\partial t}  - g \, \partial_2 z$$
	Thus, in terms of the material derivative we can write (\ref{evolutive_dam}) in the form:
\[\frac{D}{D t}(\theta \kappa)-\mbox{div}(\kappa\nabla p) =0,  \quad p\geq 0, \quad \theta \in H(p).\]
For the purpose of the time discretization, we introduce a uniform finite differences time mesh with points $t_0, t_1, \dots t_M$, with constant time step $\delta t$. Next, we introduce a forward in time approximation of the total derivative by the method of characteristics. More precisely, if we use the notation $f^n(x)=f(t^n,x)$ for $n=0, \,1, \, 2,\dots$, at each time $t^{n+1}$ we obtain the strong formulation of the discretized in time in problem,  
\begin{equation}\label{time_disc}
\frac{\theta^{(n+1)}\kappa - (\theta^{(n)}\kappa)\circ \Phi^n }{\delta t}
-\mbox{div}(\kappa\nabla p^{(n+1)}) =0,  \quad p^{(n+1)}\geq 0,  \theta^{(n+1)} \in H(p^{(n+1)})
\end{equation}
where the index $n+1$ denotes the approximation at the artificial time $t^{n+1}$ of the introduced time dependent functions and $\delta t$ denotes an artificial time step. Moreover, the function $\Phi^n$ is defined at each spatial point of the domain by $\Phi^n(x)= \Phi(t^{n+1},x; t^n)$, that denotes the position at time $t^n$ of the point placed in $x$ at time $t^{n+1}$ and moving along the integral path (characteristic curve) defined by the velocity field ${\mathbf g}$, so that $\Phi^n(x)$ can be obtained from the solution of the final value ODE problem:
$$ \displaystyle \frac{d \Phi}{d \tau}(t^{n+1},x; \tau)= {\mathbf g} \left(\tau, \Phi(t^{n+1},x; \tau)\right), \quad \Phi(t^{n+1},x; t^{n+1})=x.$$
Therefore, in terms of the previous solution we define $\Phi^n(x)=\Phi(t^{n+1},x; t^{n})$.

Taking into account the particular expression of the velocity field ${\mathbf g}$, for $x=(x_1,x_2)$ we can easily compute
\begin{equation}\nonumber 
\Phi^n(x_1,x_2)= (x_1,x_2 + g \, \Delta t ).
\end{equation}
Note that as ${\mathbf g}$ does not depend on $t$ then $\Phi^n$ does not depend on $n$, so hereafter we will drop the superindex $n$ in $\Phi^n$.

In order to write a weak form of the problem, we introduce the following functional spaces:
\begin{itemize}
\item  $V_-= H^1(D)\cap [ v|_{\Gamma_0}\leq 0]=\{\psi\in H^1(D); \psi|_{\Gamma_a}=0; \psi|_{\Gamma_0}\leq 0\}$, 
\item   $V_0= H^1_0(D,\Gamma_0)$   and note that $V_0\subset V_-$,
\item $V_+= H^1(D)\cap [ v|_{\Gamma_0}\geq0]$,
\item $V_-= H^1(D)\cap [ v|_{\Gamma_0}\leq 0]$,
\item $W_\alpha=H^1(D)\cap [v|_{\Gamma_a}=\alpha]$,
\item $W_0=H^1(D)\cap [v|_{\Gamma_a}=0]$.
\end{itemize}

Next, for $\phi \in V_{-}$ such that $\phi=0$ on $\Gamma_a$, by multiplying first equation in (\ref{time_disc}) by $\phi-p^{(n+1)}$ and integrating by parts we get the problem:\\

Find  $ p^{(n+1)} \in V_0\cap W_\alpha\cap V_+$ and $\theta^{n+1} \in L^\infty(D)$ such that 
\begin{eqnarray}
& & \displaystyle\int_D \kappa\theta^{(n+1)} (\phi-p^{n+1})
-\int_D ((\theta^{(n)}\kappa)\circ \Phi ) (\phi-p^{n+1})+ \nonumber\\ 
& & \delta t\int_D \kappa \nabla p^{(n+1)} \nabla (\phi-p^{n+1}) +\delta t\int_{\Gamma_0\cup \Gamma} 
\theta^{(n+1)}\kappa \pmb{e}_2 \pmb{n}(\phi-p^{n+1})
\geq 0, \label{ineq1}
\end{eqnarray} 
jointly with
$$ \theta^{(n+1)} \in H(p^{(n+1)}).$$

Consider now the indicatrix function of the convex set $V_{-}$ defined  over $H^1(D)$ by 
\[
I_{V_-}(v)=\left \{ \begin{array}{cc} 1 & v\in V_-,\\
+\infty & v\not \in V_-.
\end{array}
\right.
\]

Note that $I_{V_-}$ is a convex semicontinuous function so that the subdifferential operator $\partial I_{V_-}$ is a well defined maximal monotone multivalued operator, which is characterized as follows:
\begin{equation}\nonumber
    \alpha \in \partial I_{V_-}(u) \Longleftrightarrow I_{V_-}(v)-I_{V_-}(u) \geq \langle \alpha, v-u \rangle,
\end{equation}
for all $v\in H^1(D)$.

Therefore, if we define
\begin{align*}
\langle L(p^{(n+1)}), \phi-p^{(n+1)} \rangle=\displaystyle\int_D \kappa\theta^{(n+1)} (\phi-p^{n+1})
-\int_D ((\theta^{(n)}\kappa)\circ \Phi ) (\phi-p^{n+1})
+\\\delta t\int_D \kappa \nabla p^{(n+1)} \nabla (\phi-p^{n+1})+\delta t\int_{\Gamma_0\cup \Gamma}
\theta^{(n+1)}\kappa \pmb{e}_2 \pmb{n}(\phi-p^{n+1}),
\end{align*}
then from inequality (\ref{ineq1}) and the definition of the indicatrix function we get
\begin{align*}
\displaystyle
 I_{V_-}(\phi)-I_{V_-}(p^{(n+1) })
\geq  \langle L(p^{(n+1)}), \phi-p^{(n+1)} \rangle
\end{align*}
which reads   
\begin{align*}
\displaystyle\int_D \kappa\theta^{(n+1)} (\phi-p^{n+1})
-\int_D ((\theta^{(n)}\kappa)\circ \Phi ) (\phi-p^{n+1})+\delta t\int_D \kappa \nabla p^{(n+1)} \nabla (\phi-p^{n+1}) 
\\+\delta t\int_{\Gamma_0\cup \Gamma} 
\theta^{(n+1)}\kappa \pmb{e}_2 \pmb{n}(\phi-p^{n+1})
+ I_{V-}(\phi)-I_{V-}(p^{(n+1) })\geq 0
\end{align*}
for all $\phi \in W_0=H^1(D)\cap [v|_{\Gamma_a}=0]$. Moreover, the previous inequality implies that  $L(p^{(n+1)})\in \partial I_{V_-}(p^{(n+1)})$. Therefore, if we introduce the new variable $$q^{(n+1)}=L(p^{(n+1)})\in \partial I_{V_-}(p^{(n+1)})$$ 
then the time discretized problem at step $n$ can be posed as:\\

Find $  p^{(n+1)} \in X$ and  $\theta^{n+1} \in L^\infty(D)$,  such that
\begin{eqnarray}
\displaystyle\int_D \kappa\theta^{(n+1)}  \phi
+ 
\delta t\int_D \kappa \nabla p^{(n+1)} \nabla  \phi & & \nonumber\\+\delta t\int_{\Gamma_0\cup \Gamma} 
\theta^{(n+1)}\kappa \pmb{e}_2 \pmb{n} \phi 
+\delta t \int_{\Gamma_0} q^{(n+1)} \phi &=&\int_D ((\theta^{(n)}\kappa)\circ \Phi )  \phi, \label{system1}
\end{eqnarray}
for all  $\phi \in W_0$, jointly with the 
following multivalued nonlinear equations
\begin{equation}\label{nonlinearterms}
\begin{cases}
	q^{(n+1)}\in \partial I_{V_-}(p^{(n+1)}),\\ 
	\theta^{(n+1)}\in H(p^{(n+1)}).
\end{cases} 
\end{equation}

\subsection{A duality method for nonlinear terms}\label{sec:duality}
 In order to solve (\ref{system1})-(\ref{nonlinearterms}), we follow the methodology used in \cite{BermudezDurany, MR3062565} to deal with nonlinear terms associated to multivalued operators in \eqref{nonlinearterms}. These techniques are based on duality methods for nonlinear maximal monotone operators and are here applied to the multivalued Heaviside and subdifferential operators. In the seminal article \cite{BM81}, this duality methods have been introduced for solving variational inequalities. 

For this purpose, we first recall the concept of Yosida approximation. Let $G$ be a maximal monotone operator and let  $\omega$  and  $\lambda$  be non-negative real number such that $\omega\lambda<1$.  The resolvent of $G$ is defined by, 
$$J^\omega_\lambda = ( (1-\omega\lambda I)+\lambda G)^{-1}.$$
Next, we introduce the Yosida approximation of $G-\omega I $ of parameter $\lambda$, which is defined by 
\[
G_\lambda^\omega:=\frac{I-J^\omega_\lambda}{\lambda}.
\]
As it is recalled in \cite{MR3062565}, it can be proved that $\displaystyle u\in G(y)-\omega y$ is equivalent to  $u=G^\omega_\lambda (y+\lambda u)$, for further details see also the seminal article \cite{BM81}. Note that the first expression is written in terms of the multivalued operator while the second one is a nonlinear equation for $u$ in terms of an univalued Yosida operator.\\

Next, in terms of the non-negative parameters $\omega_1$ and  $\omega_2$, we introduce the new variables 
\begin{align}\label{def_alpha_beta}
&\displaystyle \alpha^{(n+1)}=q^{(n+1)}-\omega_1p^{(n+1)} 
\quad \mbox{ and }\quad
\beta^{(n+1)}=\theta^{(n+1)}-\omega_2p^{(n+1)}.
\end{align}
Therefore, from (\ref{nonlinearterms}) we have 
$$\alpha^{(n+1)}\in  \partial I_{V_-}(p^{(n+1)})-\omega_1p^{(n+1)}$$ and 
$$\beta^{(n+1)}\in  H(p^{(n+1)})-\omega_2p^{(n+1)}.$$
We can then write the variational formulation in terms of the new variables in the form
\begin{eqnarray}
\displaystyle\int_D \kappa ( \beta^{(n+1)}+\omega_2p^{(n+1)})  \phi + 
\delta t\int_D \kappa \nabla p^{(n+1)} \nabla  \phi  & & \nonumber\\+\delta t\int_{\Gamma_0\cup \Gamma} 
(( \beta^{(n+1)}+\omega_2p^{(n+1)}))\kappa \pmb{e}_2 \pmb{n} \phi & & \nonumber\\
+\delta t \int_{\Gamma_0}  (\alpha^{(n+1)}+\omega_1p^{(n+1)}) \phi 
&=&  \int_D ((\theta^{(n)}\kappa)\circ \Phi )   \phi. \nonumber
\end{eqnarray}

Next, using the previous characterization of the elements of the multivalued operator $G-\omega I$ in terms of the its Yosida approximation for the particular cases $G= \partial \partial I_{V_-}$ and $G= H$, the variational formulation can be equivalently written in the form
\begin{align}\nonumber
\int_D \kappa \nabla p^{(n+1)} \nabla \phi  
+\frac{\omega_2}{\delta t}\int_D \kappa  p^{(n+1)}  \phi+
\omega_1\int_{\Gamma_0}  p^{(n+1)} \phi+\omega_2\int_{\Gamma_0\cup \Gamma}  p^{(n+1)}\kappa \pmb{e}_2 \pmb{n} \phi \\
=
\frac{1}{\delta t}\int_D ((\theta^{(n)}\kappa)\circ \Phi )  \phi -\frac{1}{\delta t}\int_D \kappa  \beta^{(n+1)}  \phi
-\int_{\Gamma_0\cup \Gamma}  \beta^{(n+1)}\kappa \pmb{e}_2 \pmb{n} \phi 
- \int_{\Gamma_0}  \alpha^{(n+1)} \phi \label{finaleq}
\end{align}
with 
\begin{equation}
\displaystyle \alpha^{(n+1)}=(\partial I_{V_-})^{\omega_1}_{\lambda_1} (p^{(n+1)}+\lambda_1 \alpha^{n+1})\label{nl1}
\end{equation}
and 
\begin{equation}
\displaystyle \beta^{(n+1)}=H^{\omega_2}_{\lambda_2} (p^{(n+1)}+\lambda_2 \beta^{n+1}).\label{nl2}
\end{equation}

Following \cite{MR3062565}, we propose to solve numerically \eqref{finaleq}, \eqref{nl1} and \eqref{nl2} using a fixed point iteration as described in next paragraphs. 

Given  $\alpha^{(n+1)}$, $\beta^{(n+1)}$ and $\theta^{(n)}$ we solve 
equation \eqref{finaleq} for the pressure and denote the solution by 
$p^{(n+1)}=\mathcal{L}(\alpha^{(n+1)},\beta^{(n+1)},\theta^{(n)})$. So, schematically, we have the following system of coupled equations
\begin{align} \label{multi1}
\displaystyle \alpha^{(n+1)}=(\partial I_{V_-})^{\omega_1}_{\lambda_1} 
\left( \mathcal{L}(\alpha^{(n+1)},\beta^{(n+1)},\theta^{(n)})+\lambda_1 \alpha^{n+1}
\right)\\ \label{multi2}
\displaystyle \beta^{(n+1)}=H^{\omega_2}_{\lambda_2} \left(\mathcal{L}(\alpha^{(n+1)},\beta^{(n+1)},\theta^{(n)})+\lambda_2 \beta^{n+1}
\right). 
\end{align}
Using the results in \cite{MR3062565} it can be seen that, given $\theta^{(n)}$, this system can be solved by a fixed point iteration. To start the fixed point iteration we use previous values of $\alpha^{(n)}$ and $\beta^{(n)}$. The value of $\theta^n$ can be updated using \eqref{def_alpha_beta}.\\

For the spatial discretization of the linear problems arising at each step of the fixed point iteration, we consider finite elements methods. For this purpose, let $\tau^h$ be a triangulation of the domain $D$ such that it resolves the variation of the permeability coefficient $\kappa$. Consider 
$V$ the finite element space of piece-wise linear (or bi-linear) finite elements defined on the mesh $\tau$. At each step of previous iteration, the fully discretized problem can be written in terms of the solution of the following linear system: 
\begin{align}\nonumber
\left(A+\frac{\omega_2}{\delta t} M +\omega_1 M_{\Gamma_0} + 
\omega_2 M_{\Gamma_0\cup\Gamma}\right) p^{(n+1)}=\\
b^{(n)}-\left(
\frac{1}{\delta t}M+M_{\Gamma_0\cup\Gamma}\right)\beta^{(n+1)}-M_{\Gamma_0} \alpha^{(n+1)}. \label{linearproblem}
\end{align}
In the linear system (\ref{linearproblem}), we have introduced the following matrices,
$$
A=[a_{ij}] \mbox{ with } a_{ij}=\int_D \kappa \nabla \phi_i\nabla\phi_j,
$$$$
M=[m_{ij}] \mbox{ with } m_{ij}=\int_D \kappa \phi_i\phi_j,
$$$$
M_{\Gamma_0}=[m_{ij;\Gamma_0}] \mbox{ with } m_{ij;\Gamma_0}=\int_{\Gamma_0}  \phi_i\phi_j,
$$and 
$$
M_{\Gamma\cup \Gamma_0}=[m_{ij;\Gamma\cup\Gamma_0}] \mbox{ with } m_{ij;\Gamma\cup\Gamma_0}=\int_{\Gamma\cup \Gamma_0}   \phi_i\kappa \pmb{e}_2 \pmb{n} \phi_j. 
$$
Moreover, we have also considered the vector $b^{(n)} = [b_i^{(n)}]$ associated to the method of characteristics, the components of which are given by 
$$
b_i^{(n)}=
\frac{1}{\delta t}\int_D ((\theta^{(n)}\kappa)\circ \Phi )\phi_i 
$$
The computation of the integral in $b_i^{(n)}$ requires some interpolation techniques, as the function $\theta^{(n)}\kappa$ needs to be evaluated at points that may not belong to the mesh.

Note that at each step $n$ of the algorithm we recursively solve the linear system (\ref{linearproblem}) and update the terms $\alpha^{n+1}$ and $\beta^{n+1}$ in the second member by using (\ref{multi1}) and (\ref{multi2}), respectively. In practice, in the numerical examples in a forthcoming section, we consider $\lambda_1=\lambda_2=1$ and $\omega_1=\omega_2=0.5$, so that we fullfil the condition $\lambda_i \omega_i=0.5$, as in \cite{MR3062565}. Note that this condition allows to prove the convergence of the fixed point iteration in \cite{BM81} for a variational inequality problem. Also, for an elastohydrodynamic problem in magnetic storage devices the convergence is theoretically proved under the same condition in \cite{ACPV08}. We also mention that the number of fixed point iteration where chosen to be a constant number independently of of the time step iteration. This number of iteration used in the fixed point step was chosen so that the convergence to stationary solution (that is the main target of our computation) we observed. 

In next section, we propose the use of the GMsFEM to solve the fully discretized problem (\ref{linearproblem}), so that we replace the fine-scale system by a coarse linear system associated to an appropriate coarse space $V_0$.

\section{Generalized multiscale finite element method}
\label{gmsfem}
In this section we focus on high-contrast multiscale problems and summarize a GMsFEM construction 
of a coarse space $V_0$. For a more detailed description of the development of the GMsFEM methodology, see \cite{egh12,eglp13oversampling}, and references therein.

We start by choosing and initial set of basis functions that form a partition of unity. 
The space generated by this basis functions is enriched using a local spectral problem.
We use the multiscale basis functions partition of unity  with linear
boundary conditions (see \cite{eh09}, for example). We have one function per coarse-node and it is defined by 
\begin{eqnarray}\label{eq:MsFEM_standard1}
-\mbox{div}(\kappa  \nabla {\chi}_{i} )&=&0\quad \mbox{for}\ K \in \omega_i \\
 {\chi}_{i}&=& {\chi}_{i}^0  \quad 
\mbox{ on } \partial K, \nonumber
\end{eqnarray}
where ${\chi}_{i}^{0}$ is a standard linear partition of unity function.

For each coarse node neighborhood $\omega_i$, consider the eigenvalue problem
\begin{equation}
\label{eq:eig:prob}
-\mbox{div}(\kappa  \nabla \psi_\ell^{\omega_i})=\sigma_\ell^{\omega_i} \widetilde{\kappa } \psi_\ell^{\omega_i},
\end{equation}
with homogeneous Neumann
boundary condition on $\partial \omega_i$. Here $\sigma_\ell^{\omega_i}$
 and  $\psi_\ell^{\omega_i}$ are eigenvalues and
eigenvectors in $\omega_i$ and
$\widetilde{\kappa }$ is defined by
\begin{equation}\nonumber 
\widetilde{\kappa }=
\kappa \sum_{j=1}^{N_v}
{H^2}|\nabla \chi_j|^2. 
\end{equation}
We use an ascending ordering on the eigenvectors, 
$
\sigma_1^{\omega_i}\leq \sigma_2^{\omega_i}\leq....
$\\

Using the partition of unity functions from Eq. \eqref{eq:MsFEM_standard1} and eigenfunctions from Eq. \eqref{eq:eig:prob}, we then construct a set of enriched multiscale basis functions given by
$\chi_i \psi^{\omega_i}_\ell$ for selected eigenvectors $\psi^{\omega_i}_\ell$. Using
$L_i$ to denote the number of basis functions from the coarse region $\omega_i$, we then define the coarse GMsFEM space by
\[V_0=\mbox{span}\{\Phi_{i,\ell}=\chi_i \psi_\ell^{\omega_i},\quad i=1,\dots,N_v,\quad \ell=1,\dots,L_i\}. \]
For more details, motivation of the construction, and approximation properties of the space $V_0$ as well as the choice of the initial partition of unity basis functions we refer the interested reader to \cite{egh12}.

Summarizing, in order to solve problem (\ref{linearproblem}) for the pressure we use the GMsFEM coarse space $V_0$ constructed in 
this section. More precisely, let $R_0$ the matrix whose columns correspond to the coarse basis functions, that is, the column space of $R_0$ if $V_0$. Instead of solving the fine-scale linear system \eqref{linearproblem} we solve the coarse-linear system 
\begin{equation}\label{eq:coarse_system}
S_0 p_0^{(n+1)} = c_0^{(n)} \,,
\end{equation}
where the matrix $S_0$ and the second member $c_0^{(n)}$ are given by
\begin{eqnarray}
S_0 & = & R_0^T\left(A+\frac{\omega_2}{\delta t} M +\omega_1 M_{\Gamma_0} + \omega_2 M_{\Gamma_0\cup\Gamma}\right) R_0 \nonumber\\  
c_0^{(n)} & = & R_0^T \left(b^{(n)}-\left(
\frac{1}{\delta t}M+M_{\Gamma_0\cup\Gamma}\right)\beta^{(n+1)}-M_{\Gamma_0} \alpha^{(n+1)}\right).\nonumber 
\end{eqnarray} 
We then maintain the duality method explained before but using the approximation $R_0p_0^{(n+1)}\approx p^{(n+1)}$ that makes the computation more efficient since, instead of solving the full resolution linear system 
\eqref{linearproblem}, we solve the small coarse problem 
\eqref{eq:coarse_system}. We mention that due to the high-contrast multiscale structure of the coefficient, we need to solve the coarse problem at the right resolution in order to obtain good approximation; see \cite{egh12}. With the GMsFEM methodology we can adapt the resolution of the coarse solver in order to be able to obtain good results with the duality method. 

We recall that, for each time iteration and each fixed point iteration, in order to compute the current pressure a linear system has to be solved. In the here proposed methodology, instead of solving the fine-grid linear system \eqref{linearproblem} we solve the coarse scale linear system \eqref{eq:coarse_system}. 
We stress that linear system \eqref{linearproblem} is very large and ill-conditioned (with condition number increasing with the contrast in the coefficient).
The size of system 
\eqref{eq:coarse_system} is of the order of the number of coarse-scale nodes so it is suitable for factorization methods. This allows us to save computational time. Moreover, we mention that not the basis functions nor the coarse scale operators and matrices change  throughout
the time and non-linear iteration. Therefore, the set-up cost (constructing coarse grid, computing local eigenvectors and coarse basis functions and assembling coarse scale operators) can be consider as a pre-processing cost. 
See \cite{egh12} for more details on the computational implementation of GMsFEMs. 

A key aspect of the GEMsFEM is that the needed resolution  can be considered a priori or a posteriori depending on the application. 
In this paper we show how the resolution of the method, that is, the parameters $L_i$ determining the number of eigenvectors in $\omega_i$ used in the construction of the coarse space, affects the solution of the free boundary of the heterogeneous multiscale dam problem. As in others applications leading to iterative corrections using solutions of the diffusion equations, if the approximation of this step is poor, overall poor results are obtained in the solution 
of the problems.

In the numerical results of the next section we observed that the  number of time steps needed to obtain convergence to stationary solution (using the same tolerance of the increment) is larger when using the GMsFEM approximation. It was needed around 10\% more time steps that in the case of computing the fine scale solution \eqref{linearproblem}. Recall that when using the GMsFEM approximation we solve 
\eqref{eq:coarse_system} that is is very small dimension when compared to the fine-scale system. See Table 
\ref{taberrors} below. 

\section{Numerical results}\label{sec:numerics}
In this section we present numerical illustrations to show the performance of the GMSFEM method when compared to the reference solution. We recall that the time evolution and the duality method are maintained as before using vectors in the fine grid but instead of computing fine-scale solves for the pressure equations we use the GMsFEM solution. 
We consider the following boundary partition (see Figure 
\ref{illustration})
\begin{align*}
\Gamma_0&=\{0\}\times [3/5,1]\cup [0,1]\times\{1\}
\cup \{1\}\times [2/5,1],\\
\Gamma_a&=\{0\}\times [0,3/5]\cup 
 \{1\}\times [0,2/5],\\
 \Gamma&=[0,1]\times  \{0\},
\end{align*}
and we have the Dirichlet data $p=4/5-x_2$  on $\{0\}\times [(3/5),1]$ and
 $p=1/5-x_2$  on $\{1\}\times [2/5,1]$ and $p=0$ 
 on $\Gamma_0$.
 
\begin{figure}[h]
\includegraphics[width=.4\textwidth]{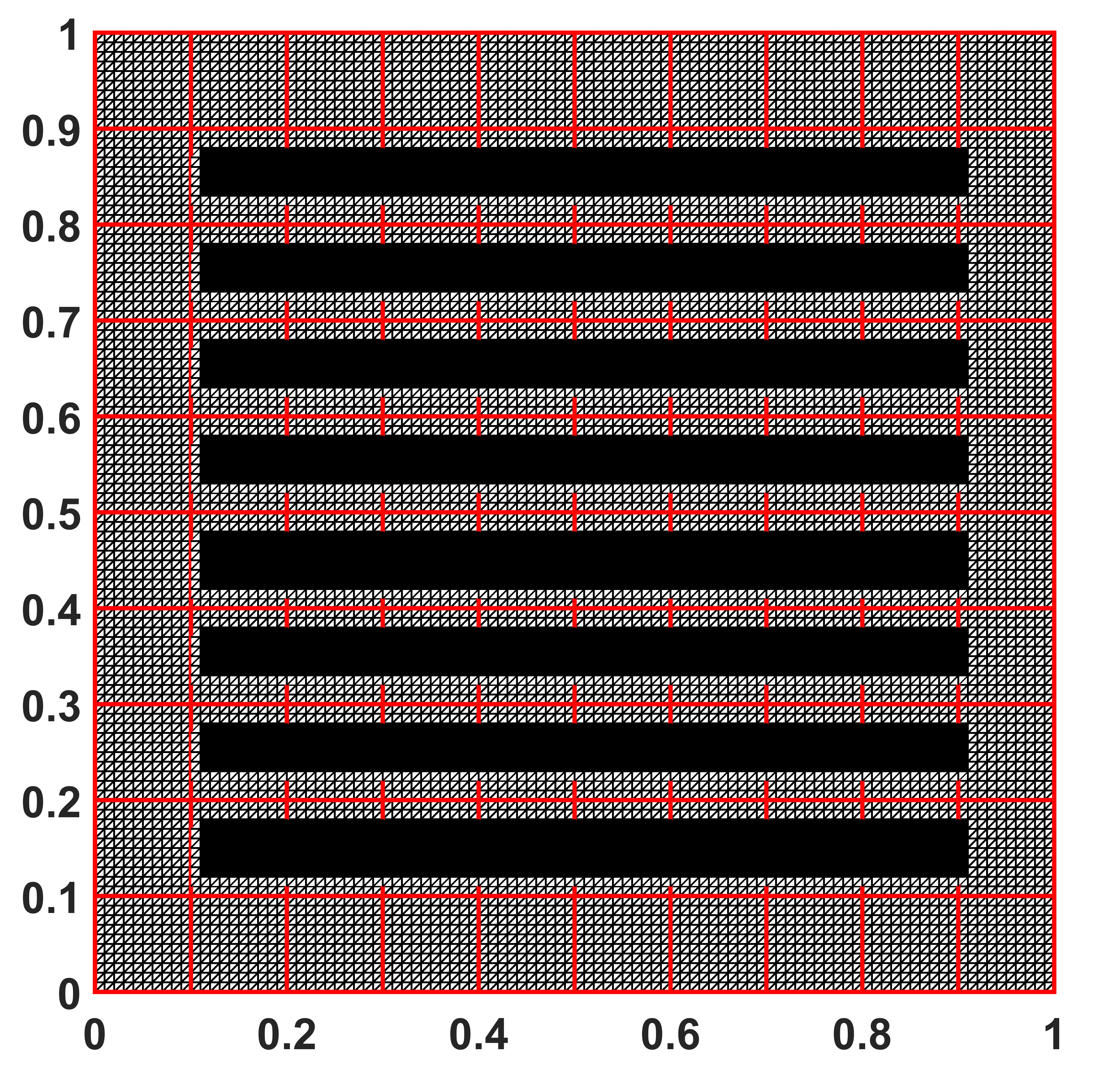}
\includegraphics[width=.4\textwidth]{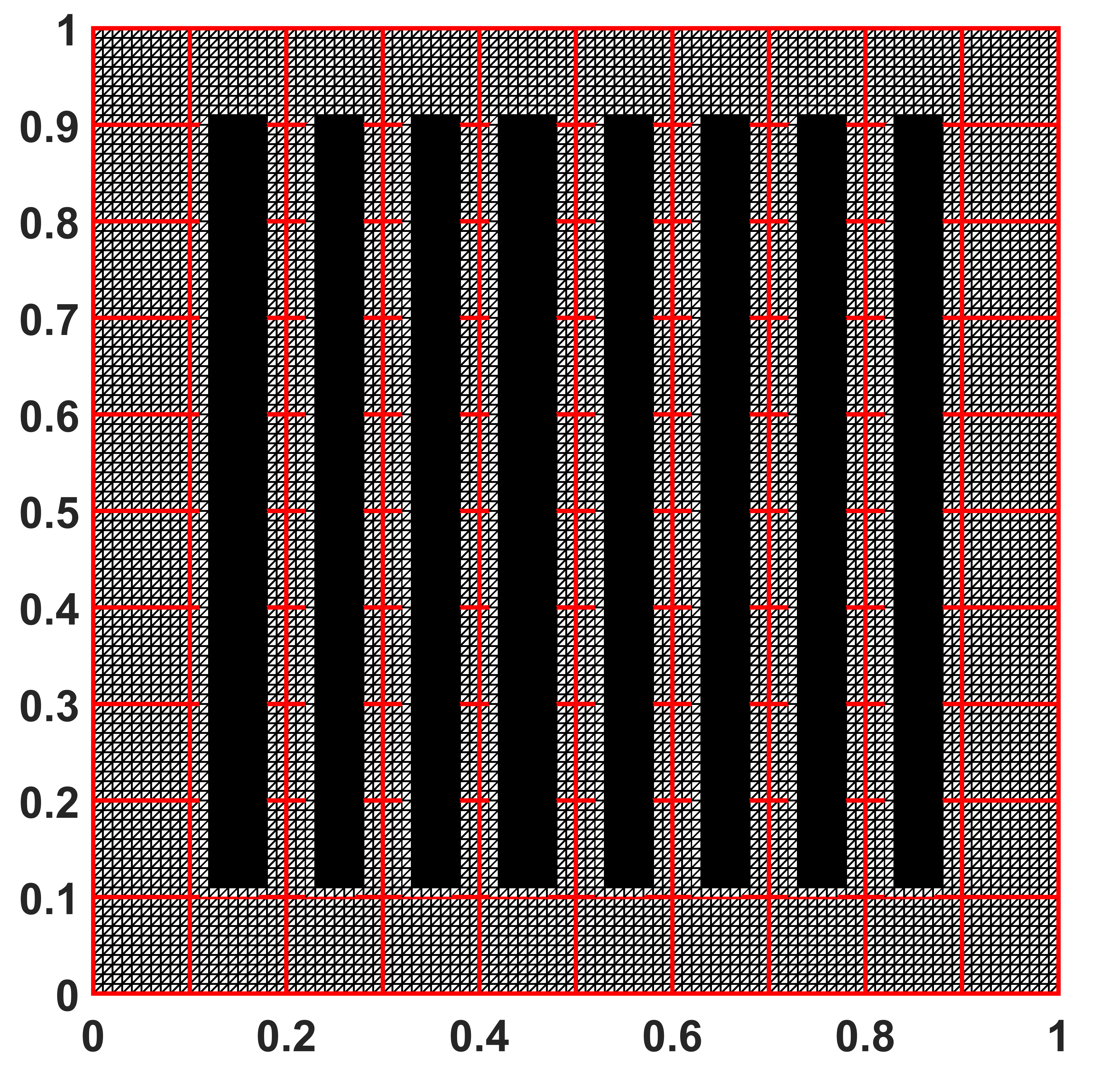}\\
\includegraphics[width=.4\textwidth]{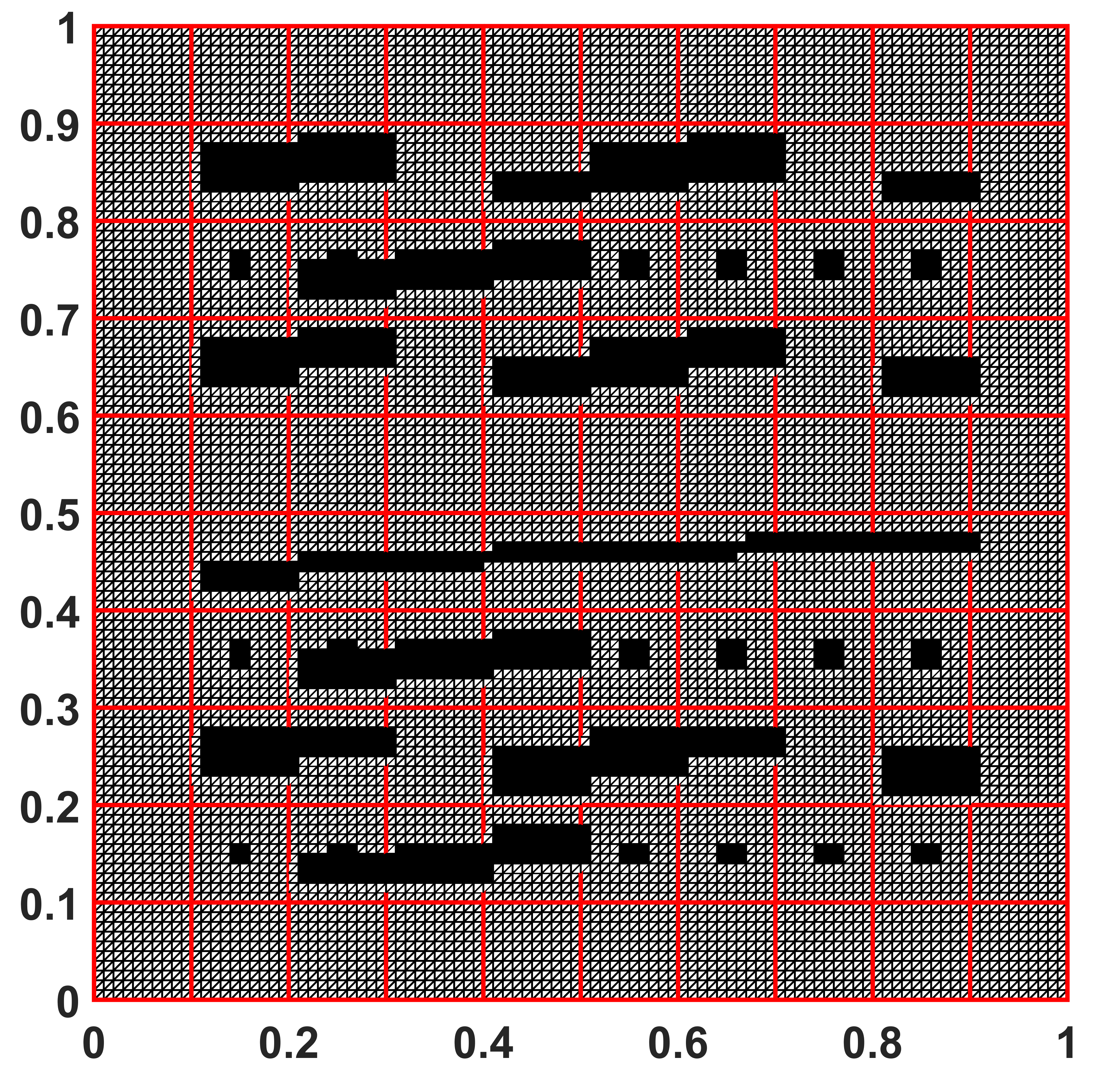}\\
\caption{ High-Contrast coefficients used in the numerical experiment. 
High-conductivity channels in solid black. Top-Left picture: coefficient with horizontal channels and non-aligned coarse mesh. Top-Right figure: coefficient with vertical channels and non-aligned coarse grid. Bottom figure: A high-contrast coefficient with channels and inclusions. In our numerical experiments we use coefficients of background 1 and high-contrast value  $10^2$.\label{fig:highcontrast}}
\end{figure}

In order to numerically study the performance of the GMsFEM method applied to the heterogeneous dam problems we use the coefficients 
depicted in Figure \ref{fig:highcontrast}. We then compare to the reference solutions, that is, we compute the error between the solution of the overall iteration with solutions of the diffusion equation on the fine-gird, with the multiscale solutions, that is, the solutions obtained by using the coarse-scale solution $p_0^{(n+1)}$ (downscaled to the fine-grid as $R_0p_0^{(n+1)}$) for the approximation of the diffusion equation. 
In particular we use $L_i=0,1,\dots,10$ for all $i$. 
We run the time iteration and the fixed point iteration until the norm of the increment is less than a given tolerance ($10^{-4}$ in our numerical test).

We consider a structured fine-grid with 100 elements in each direction (yielding a fine-scale linear system matrix of 
dimension 10000$\times$10000). We also consider a coarse mesh (made of squares) with 10 elements in each direction. The coarse mesh is a structured mesh and it is not aligned to  the variations of the coefficient, see 
\cite{EfendievGLWESAIM12,eglp13oversampling} and 
Section \ref{sec:2.1}. We consider three different piece-wise constant coefficients:
horizontal channels, vertical channels and a high-contrast multiscale coefficient with channels and inclusions. In our numerical experiments we use coefficients of background 1 and high-contrast value  $10^2$.
See Figure \ref{fig:highcontrast}.

\subsection{Horizontal channels}
\begin{figure}[h]
\centering
\includegraphics[width=.4\textwidth]{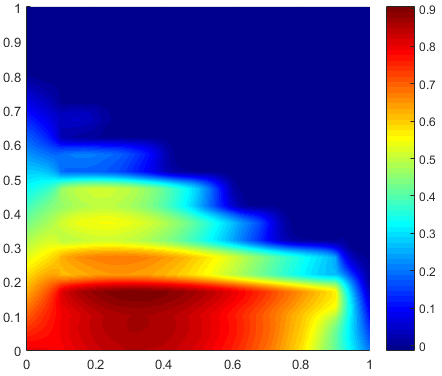}
\includegraphics[width=.4\textwidth]{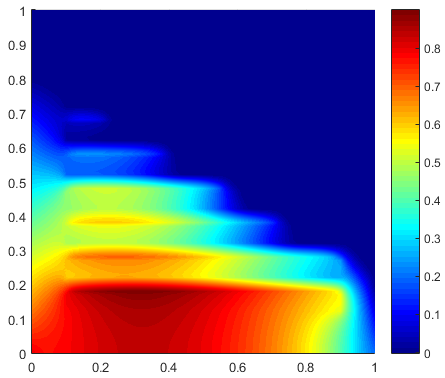}\\
\includegraphics[width=.4\textwidth]{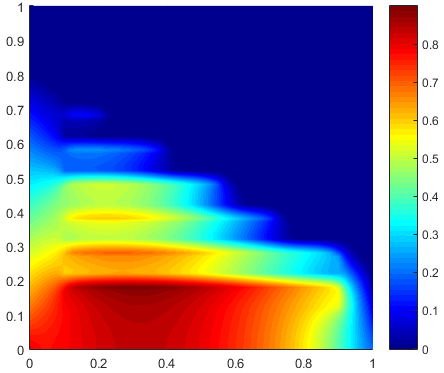}
\includegraphics[width=.4\textwidth]{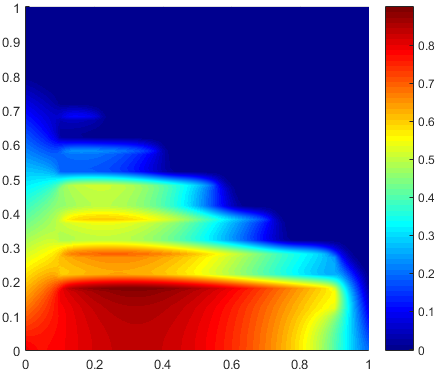}
\caption{Computed pressure for horizontal channels coefficient: Fine-grid solution  and coarse scale solutions for a variety of coarse space dimensions. From top to bottom and left to right: Using fine-grid solution, $dim(A)=101^2\times 101^2$. Using coarse-grid solution with $L_i=1$, $dim(A_0)=11^2\times 11^2 $.
Using coarse-grid solution with $L_i=2$, $dim(A_0)=(2*11^2)\times
(2*11^2) $.
Using coarse-grid solution with $L_i=4$, $dim(A_0)=(4*11^2)\times
(4*11^2)$.
\label{diffhorizontal}}
\end{figure}

\begin{figure}[h]
\includegraphics[width=.4\textwidth]{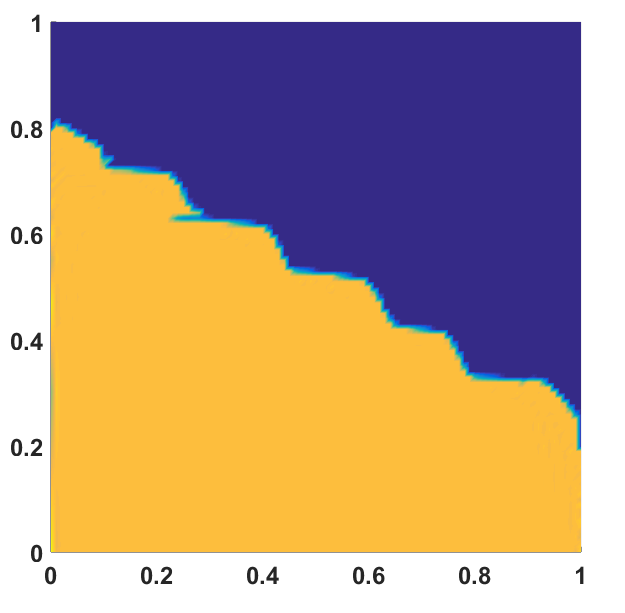}
\includegraphics[width=.4\textwidth]{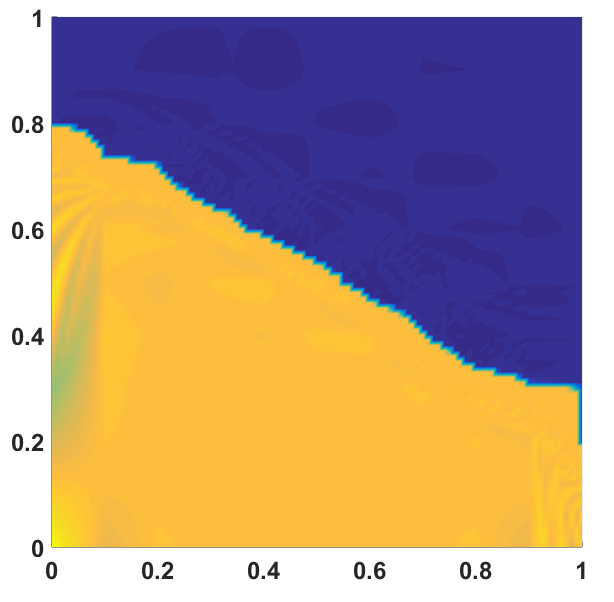}\\
\includegraphics[width=.4\textwidth]{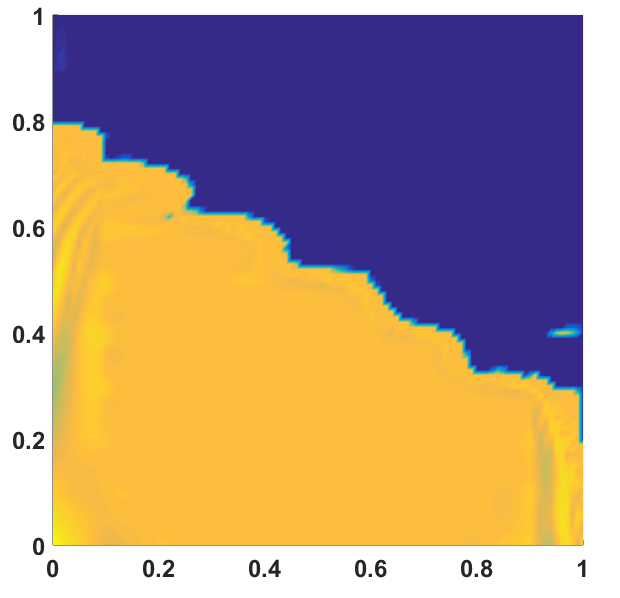}
\includegraphics[width=.4\textwidth]{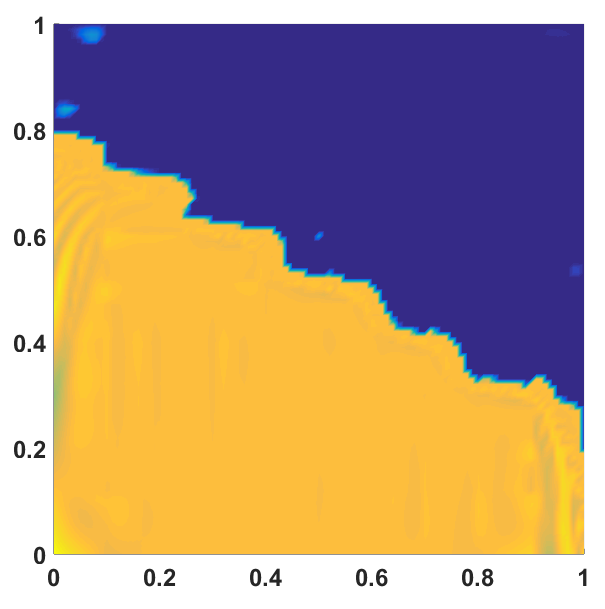}
\caption{Computed saturation for horizontal channels coefficient: Fine-grid solution (Top-Left figure) and coarse scale solutions for a variety of coarse space dimensions. From top to bottom and left to right: Using fine-grid solution, $dim(A)=101^2\times 101^2$. Using coarse-grid solution with $L_i=1$, $dim(A_0)=11^2\times 11^2 $.
Using coarse-grid solution with $L_i=2$, $dim(A_0)=(2x*11^2)\times
(2*11^2) $.
Using coarse-grid solution with $L_i=4$, $dim(A_0)=(4*11^2)\times
(4*11^2)$. .\label{sathorizontal}}
\end{figure}

In our first numerical computations  we consider the horizontal channels case; see 
Figure \ref{fig:highcontrast}. 
See Figure \ref{diffhorizontal} for the pressure solution and   Figure \ref{sathorizontal} for the computed saturation for a variety of coarse 
spaces dimensions. 

\subsection{Vertical channels}
\begin{figure}[h]
\centering
\includegraphics[width=.4\textwidth]{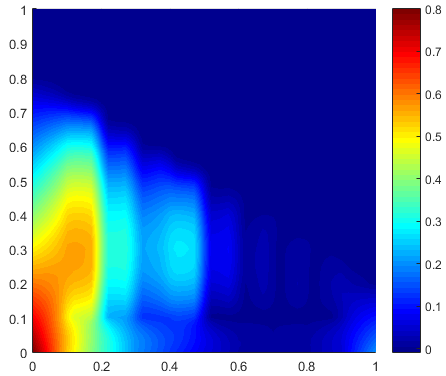}
\includegraphics[width=.4\textwidth]{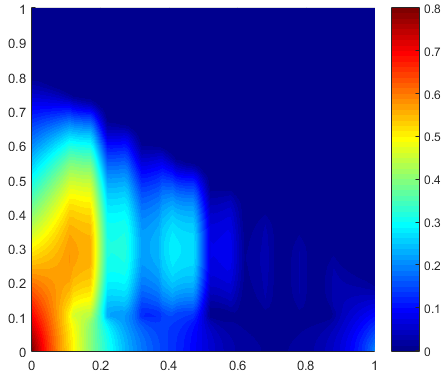}\\
\includegraphics[width=.4\textwidth]{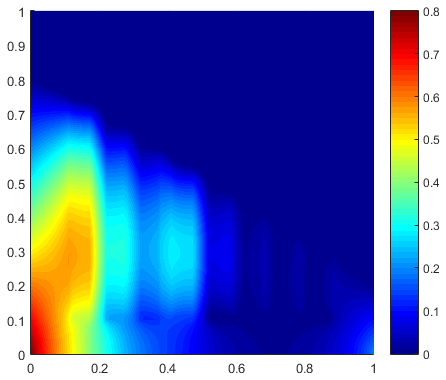}
\includegraphics[width=.4\textwidth]{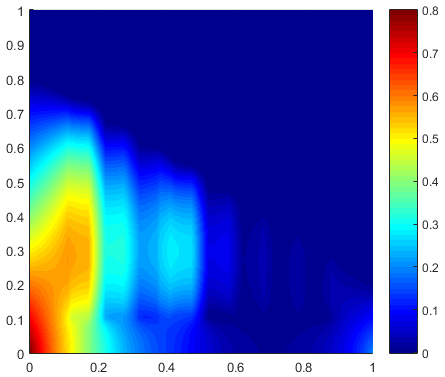}
\caption{Computed pressure for vertical channels coefficient: Fine-grid solution (Top-Left figure) and coarse scale solutions for a variety of coarse space dimensions. From top to bottom and left to right: Using fine-grid solution, $dim(A)=101^2\times 101^2$. Using coarse-grid solution with $L_i=1$, $dim(A_0)=11^2 \times 11^2 $.
Using coarse-grid solution with $L_i=2$, $dim(A_0)=(2*11^2)\times
(2*11^2) $.
Using coarse-grid solution with $L_i=4$, $dim(A_0)=(4*11^2)\times
(4*11^2)$..\label{difvertical}}
\end{figure}

\begin{figure}[h]
\includegraphics[width=.4\textwidth]{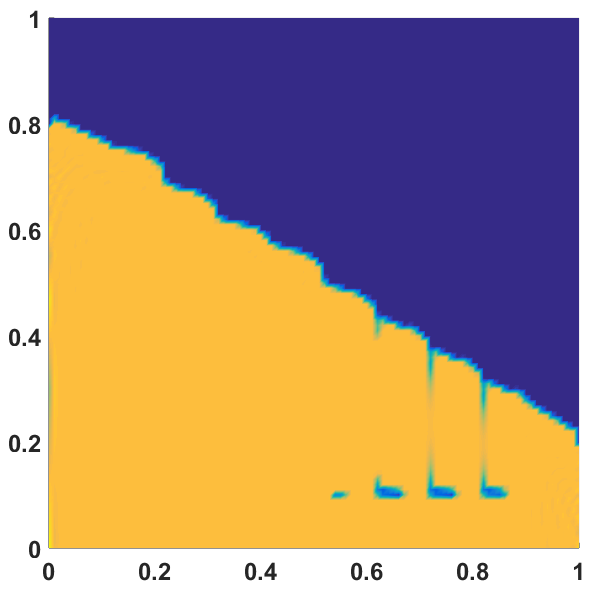}
\includegraphics[width=.4\textwidth]{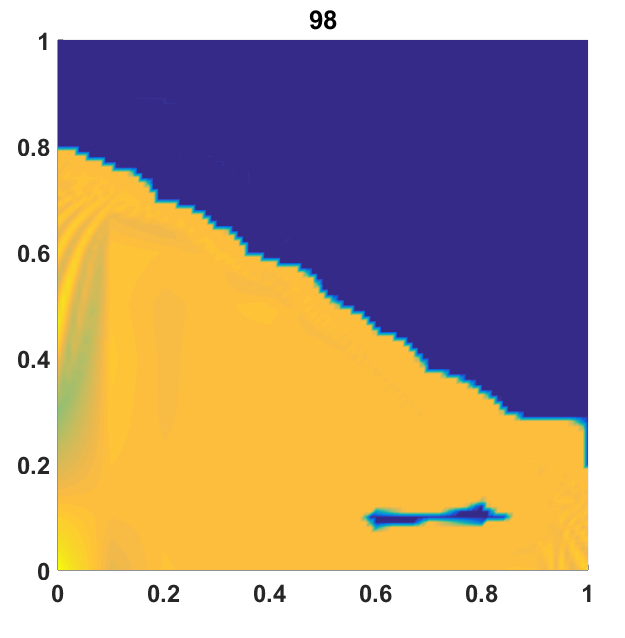}\\
\includegraphics[width=.4\textwidth]{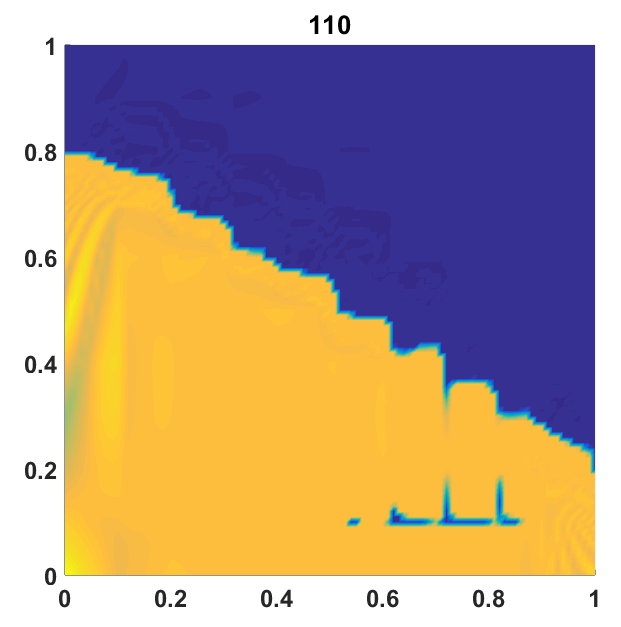}
\includegraphics[width=.4\textwidth]{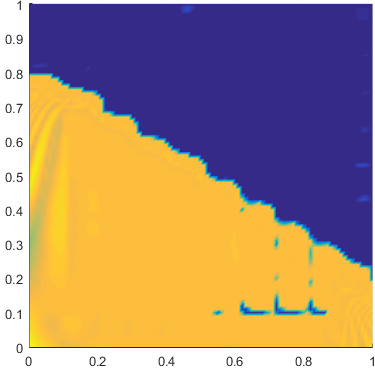}
\caption{Computed solutions for vertical channels coefficient: Fine-grid solution (Top-Left figure) and coarse scale solutions for a variety of coarse space dimensions. From top to bottom and left to right: Using fine-grid solution, $dim(A)=101^2\times 101^2$. Using coarse-grid solution with $L_i=1$, $dim(A_0)=11^2\times 11^2 $.
Using coarse-grid solution with $L_i=2$, $dim(A_0)=(2*11^2)\times
(2*11^2) $.
Using coarse-grid solution with $L_i=4$, $dim(A_0)=(4*11^2)\times
(4*11^2)$..\label{satvertical}}
\end{figure}

We consider now  a vertical channels case;
see Figure \ref{fig:highcontrast}. We obtain the results in Figure \ref{difvertical} for the pressure solution and the results in Figure \ref{satvertical} for the computed saturation for different coarse space 
dimensions. 

\subsection{A high-contrast coefficient with channels and inclusions}
\begin{figure}[h]
\centering
\includegraphics[scale=.53]{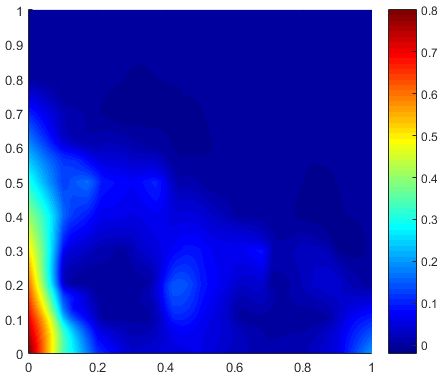}
\includegraphics[scale=.53]{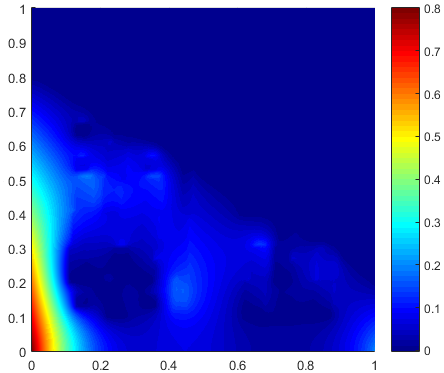}\\
\includegraphics[scale=.53]{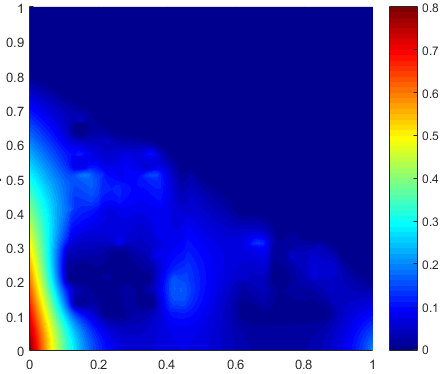}
\includegraphics[scale=.53]{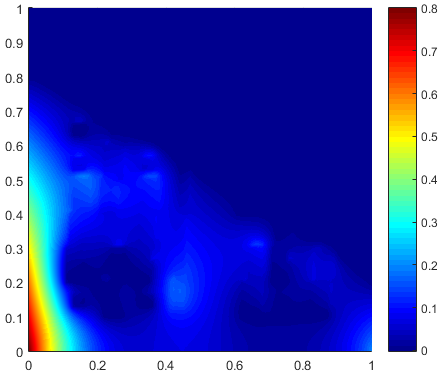}
\caption{Computed pressure for coefficient with channels and inclusions: Fine-grid solution (Top-Left figure) and coarse scale solutions for a variety of coarse space dimensions. From top to bottom and left to right: Using fine-grid solution, $dim(A)=101^2\times 101^2$. Using coarse-grid solution with $L_i=1$, $dim(A_0)=11^2 \times 11^2 $.
Using coarse-grid solution with $L_i=2$, $dim(A_0)=(2*11^2)\times
(2*11^2) $.
Using coarse-grid solution with $L_i=4$, $dim(A_0)=(4*11^2)\times
(4*11^2)$..\label{diffinclusions}}
\end{figure}

\begin{figure}[h]
\centering
\includegraphics[scale=.56]{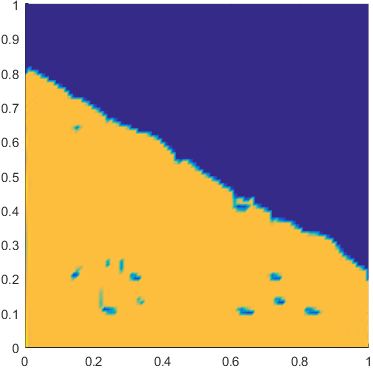}
\includegraphics[scale=.56]{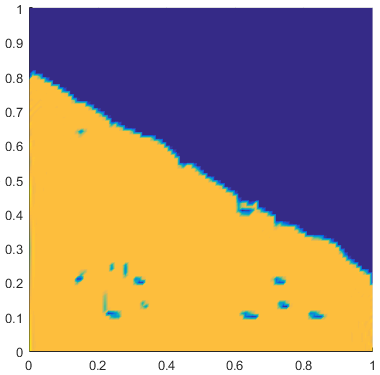}\\
\includegraphics[scale=.5]{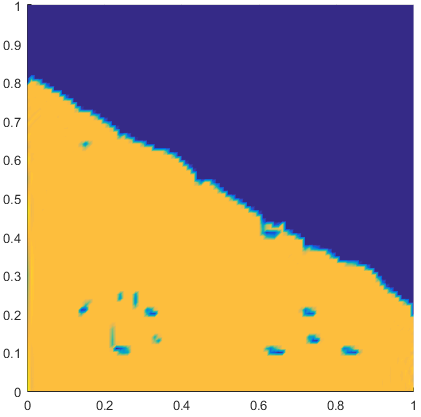}
\includegraphics[scale=.5]{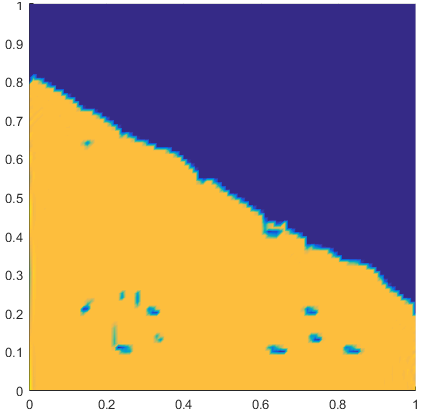}
\caption{Computed solutions for coefficient with channels and inclusions: Fine-grid solution (Top-Left figure) and coarse scale solutions for a variety of coarse space dimensions. From top to bottom and left to right: Using fine-grid solution, $dim(A)=101^2\times 101^2$. Using coarse-grid solution with $L_i=1$, $dim(A_0)=11^2\times 11^2 $.
Using coarse-grid solution with $L_i=2$, $dim(A_0)=(2*11^2)\times
(2*11^2) $.
Using coarse-grid solution with $L_i=4$, $dim(A_0)=(4*11^2)\times
(4*11^2)$..\label{satinclusions}}
\end{figure}

We consider the high-contrast coefficient with channels and inclusions depicted 
in Figure \ref{fig:highcontrast}. We obtain the results in Figure \ref{diffinclusions} for the pressure solution and 
the results in Figure \ref{satinclusions} for the computed saturation for different 
coarse space dimensions. For this examples we show the percentage of error (in the energy norm) when comparing the fine-scale solution and the GMsFEm solution. 
See Table \ref{taberrors}. We see that the approximation improves as we add more basis functions, that is, as we use more eigenvectors in the construction of the coarse problem.

\begin{table}
  \centering
  \begin{tabular}{  | c | c  |}
    \hline
  Coarse Dimension  & Relative Errors (\%) \\
    \hline \hline
      $\text{dim}(V_0)$  [$\#$ basis]  & $H^1_k (\Omega)$ \\
    \hline
   121  [1]   &    16.31  \\
    \hline
     202  [2]   &        13.51     \\
    \hline
     364  [4]  &   11.59         \\
    \hline
      526  [6]  & 10.23           \\
    \hline
     688  [8]  &   9.63        \\
    \hline
     850  [10]  & 8.76      \\
    \hline
  \end{tabular}
\caption{Relative pressure errors for a variety of coarse space dimensions for coefficient with channels and inclusions.}
\label{taberrors}
\end{table}

\end{document}